\documentclass[12pt]{article}
\usepackage{amsmath,amssymb,amsfonts}
\usepackage{graphicx}
\newtheorem{thm}{Theorem}
\newtheorem{lem}{Lemma}

\title{ H-Infinity Optimal Decentralized Matching Model\\
Is Not Always Rational}

\author{Alexandre Megretski
\thanks{
Department of Electrical Engineering and Computer Science at the Massachusetts Institute of Technology, Cambridge, MA 02149. email: ameg@mit.edu.} 
}

\begin{document}

\maketitle

\begin{abstract}
We construct structured H-Infinity optimal model matching
 problems with rational coefficients, in which the optimal 
solution is not rational, in the sense that
the cost does not achieve its  maximal lower bound 
on the set of rational 
matching models, but the same infimum can be reached by
using a continuous non-rational matching model.
\end{abstract}

\subsection*{Notation and terminology}%
We use $\mathbb C$ to denote the set of all complex numbers, and
$\mathbb C_{m\times n}$ to denote the set of all $m$-by-$k$
complex matrices. For $M\in\mathbb C_{m\times n}$,
the element in the $i$-th row and the $r$-th
column of $M$ is expressed by $[M]_{i,r}$. 
For a positive integer $n$, $I_n$ is the $n$-by-$n$ identity matrix, and
$J_n$ denotes the $n$-by-$n$ "order reversal" matrix (i.e.
$[J_n]_{i,r}=1$ when $i+r=n+1$, and $[J_n]_{i,r}=0$ 
otherwise).
We use
\[
\mathbb D=\{w\in\mathbb C:~|w|<1\},\quad
\mathbb T=\{w\in\mathbb C:~ |w|=1\},\quad
\mathbb D_+=\{w\in\mathbb C:~ |w|>1\}
\]
to denote the open unit disc, its boundary, and the complement
of its closure.
For all $w\in\mathbb C$, $\bar w$ is the complex
conjugate of $w$, and, 
for $w=re^{j\theta}$, where  $-\pi<\theta\le\pi$, and
$\beta>0$, the power $w^\beta$ is defined by 
$w^\beta=r^\beta e^{j\beta\theta}$. As a shortcut, we use
$\sqrt{w}$ for $w^{1/2}$.

When $\Omega$ is a subset of $\mathbb C$, a function
$G:~\Omega\to\mathbb C$ is said to be {\sl real symmetric}
when $\bar w\in\Omega$ and $G(\bar w)=\overline{G(w)}$
whenever $w\in\Omega$.
$\mathbf H^\infty$ denotes 
the set of  real symmetric analytic functions
$G:~\mathbb D\to\mathbb C$ such that
$\sup_{w\in\mathbb D}|G(w)|<\infty$,
while
$\mathbf A$ is the subset of functions $G\in\mathbf H^\infty$
which can be extended continuously to $\mathbb D\cup\mathbb T$,
and $\mathbf{RA}$ is the subset of {\sl rational}
functions $G\in\mathbf A$.
Furthermore, $\mathbf H^\infty_{m\times k}$, 
$\mathbf A_{m\times k}$, and
$\mathbf{RA}_{m\times k}$ denote the sets of
$m$-by-$k$  complex matrix-valued functions
$G:~\mathbb D\to\mathbb C_{m\times n}$ for which all
components $G_{i,r}:~z\to[G(z)]_{i,r}$ belong to the classes 
$\mathbf H^\infty$, $\mathbf A$, and
$\mathbf{RA}$, respectively. 
The elements of $\mathbf H^\infty_{m\times k}$ will be referred to as
{\sl stable transfer matrices} (or {\sl stable transfer functions}
in the case $m=k=1$), though it is more common in control systems
literature to call so the functions $F:~D_+\to\mathbb C$
defined by $F(z)=G(1/z)$ for some 
$G\in\mathbf H^\infty_{m\times k}$.
Naturally, transfer matrices from 
$\mathbf A_{m\times k}$ will be viewed as 
continuous functions on $\mathbb D\cup\mathbb T$.
For 
$G\in\mathbf H^\infty_{m\times k}$, the {\sl L-Infinity norm}
$\|G\|_{\infty}$ is defined as  the minimal upper bound of the 
largest singular number $\sigma_{\max}(G(w))$
of $G(w)$ over $w\in\mathbb D$.
Every function $G\in\mathbf{RA}_{m\times k}$ can be represented 
(in many ways) in the form $G(w)=D+wC(I_n-wA)^{-1}B$, where
$A,B,C,D$ are real matrices of dimensions $n$-by-$n$, $n$-by-$k$,
$m$-by-$n$, and $m$-by-$k$ respectively (when $n=0$, the
representation is interpreted as $G(w)\equiv D$). The minimal possible
value of $n$ in such representation will be referred to as the {\sl order}
of $G$, and the set of all elements $G\in\mathbf{RA}_{m\times k}$
of order not larger than $n$ will be denoted by 
$\mathbf{RA}^{n}_{m\times k}$.

\section{Introduction}
Given rational stable transfer matrices 
$L_0\in\mathbf{RA}_{p\times d}$, 
$L_1\in\mathbf{RA}_{p\times m}$,  
$L_2\in\mathbf{RA}_{k\times d}$, 
the classical {\sl  H-Infinity model matching problem} can be expressed
in the form
\begin{equation}\label{eq0}
\|L_0+L_1QL_2\|_\infty\to\min_{Q\in\mathbf H^{\infty}_{m\times k}}.
\end{equation}
In other words, it
calls for finding
a stable transfer matrix 
$Q\in\mathbf H^\infty_{m\times k}$ which minimizes 
{\sl model matching error} $\|L_0+L_1QL_2\|_\infty$. 
In this paper, we only consider the case
when the {\sl well-posedness} assumption
\begin{equation}\label{eq1}
L_1(z)'L_1(z)>0,\quad L_2(z)L_2(z)'>0\qquad\forall~z\in\mathbb T,
\end{equation}
is satisfied, thus guaranteeing existence of an optimal 
$Q\in\mathbf H^\infty_{m\times k}$.

The classical H-Infinity model matching problem
is well studied, as it
appears naturally (after applying the
so-called "Youla-", or "Q-", parameterization)
as an intermediate step in designing a stabilizing
linear time invariant feedback  for a given stabilizable
finite order linear time invariant plant (with $d$ noise
inputs,  $m$ actuator inputs, $p$ cost outputs, and
$k$ sensor outputs), with an objective of minimizing the L2
induced norm in the closed loop map 
from the noise inputs to the cost outputs.
In particular, when
\[   L=\begin{bmatrix} L_0 & L_1\\ L_2 & 0\end{bmatrix}\in
\mathbf{RA}^n_{(p+k)\times(d+m)}\]
has order $n$,
restricting $Q$ to be rational of order 
$n$ does not reduce the best achievable
performance, in the sense that 
\[ \min_{Q\in\mathbf{RA}^n_{m\times k}}\|L_0+L_1QL_2\|_\infty=
\min_{Q\in\mathbf H^{\infty}_{m\times k}}\|L_0+L_1QL_2\|_\infty.\]

In the last decade, breakthrough advances in understanding
Q-parameterization (see, for example, \cite{Lall06}) led naturally to
a {\sl structured} version of the model matching problem
(\ref{eq0}), in which some entries of $Q$ are constrained to be 
identically zero.
Such formulations are obtained, for plants of a special structure,
when there is a need to optimize a
{\sl decentralized} stabilizing linear time invariant
feedback.

One basic question associated with this development is whether a
{\sl rational} optimal $Q$ is guaranteed 
 to exist (subject to assumption (\ref{eq1})) in the problem
of minimizing the cost $\|L_0+L_1QL_2\|_\infty$ when
$Q$ is restricted to the set 
of all {\sl diagonal}
stable transfer matrices of appropriate dimension. 
This paper aims to answer the question 
(posed to the author by S. Lall) negatively. 

Specifically, let 
$\mathcal D$ denote the set of all {\sl diagonal} stable transfer
matrices. We produce triplets
$(L_0,L_1,L_2)$ of rational stable transfer matrices 
$L_i\in\mathbf{RA}_{2\times 2}$  satisfying condition
(\ref{eq1}), such that
\[ \inf_{Q\in\mathbf{RA}^{2\times 2}\cap\mathcal D}
\|L_0+L_1QL_2\|_\infty=
\min_{Q\in\mathbf H^\infty_{2\times 2}\cap\mathcal D}
\|L_0+L_1QL_2\|_\infty\]
is achieved at a unique 
$Q\in\mathbf A_{2\times 2}\cap\mathcal D$ which is
 is {\sl not} a rational function (in fact, the optimal
$Q$ can be computed explicitly).
The derivation relies on the conformal mapping technique
by Allen Tannenbaum \cite{Tannenbaum80}.

\section{Main Results}
We will use the standard expression for
the conformal map of the
open unit disc $\mathbb D$ to the 
open "lens" region
\begin{equation}
 \Omega_\gamma=\left\{s\in\mathbb C:\quad |1-s|<
\gamma,\quad
|1+s|<\gamma\right\}\qquad\qquad(1<\gamma<\infty).
\end{equation}

\begin{lem}\label{lem1}
For every $\alpha\in(0,\pi/2)$ and 
$\gamma=\frac1{\cos\alpha}\in(1,\infty)$ the function $F_\gamma:~\mathbb D\to\mathbb C$
defined by
\[F_\gamma(w)=
j\tan(\alpha)\cdot
\frac{1-\left(\frac{1+jw}{1-jw}\right)^{2\alpha/\pi}}{1+
\left(\frac{1+jw}{1-jw}\right)^{2\alpha/\pi}}
\qquad\qquad (w\in\mathbb D)
\]
belongs to class $\mathbf A$,  establishes
a bijection between $\mathbb D$ and $\Omega_\gamma$,
 and satisfies the condition
$\dot F_\gamma(0)=2\alpha\tan(\alpha)/\pi$.
\end{lem}

The following statement provides a simple example of a  structured
H-Infinity optimal model matching problem with $d=p=m=k=2$, 
$L_2=I_2$, such that
the optimal $Q\in\mathbf H^\infty_{2,2}\cap\mathcal D$ 
is unique, belongs to the class $\mathbf A^\infty_{2,2}$,
but is not a rational function.

\begin{thm}\label{thm0}
Equalities
\[ \inf_{Q\in\mathbf{RA}_{2\times2}\cap\mathcal D}
\|L_0+L_1Q\|_\infty=
\min_{Q\in\mathbf H^\infty_{2\times2}\cap\mathcal D}
\|L_0+L_1Q\|_\infty=\sqrt{2}\]
hold for 
\[  L_0(w)=I_2+0.5wJ_2=
\begin{bmatrix} 1& 0.5w\\ 0.5w&1\end{bmatrix},\quad
L_1(w)=w^2J_2=
\begin{bmatrix}  0& w^2\\ w^2 & 0\end{bmatrix}.\]
Moreover,
the only $Q_*\in\mathbf H^\infty_{2\times2}\cap\mathcal D$
such that $\|L_0+L_1Q_*\|_\infty=\sqrt2$ is given by
$Q_*(w)=S_*(w)I_2$, where $S_*\in\mathbf{A}$ is defined
by $0.5w+w^2 S_*(w)=F_{\sqrt{2}}(w)$, and
$F_\gamma\in\mathbf A$ is defined in Lemma~\ref{lem1}.
\end{thm}

A proof of Theorem~\ref{thm0} is given in the Appendix section below.

The optimization task described in Theorem~\ref{thm0}
is actually a special case of a slightly more general class of structured
model matching problems in which the optimal
diagonal $Q$ is guaranteed to be continuous but not rational.

\begin{thm}\label{thm1}
Let $a,b\in\mathbf{RA}$ be 
such that $b(w)\neq0$ for all $w\in\mathbb T$, and
$a+bq$ is not constant for every $q\in\mathbf{RA}$.
Then, for
$L_0(w)=I_2+a(w)J_2$, 
$L_1(w)=b(w)J_2$, 
the only
$Q_*\in\mathbf H^\infty_{2\times2}\cap\mathcal D$
satisfying 
\[  \|L_0+L_1Q_*\|_\infty=
\inf_{Q\in\mathbf{RA}_{2\times2}\cap\mathcal D}
\|L_0+L_1Q\|_\infty=
\min_{Q\in\mathbf H^\infty_{2\times2}\cap\mathcal D}
\|L_0+L_1Q\|_\infty=\gamma>1\]
is given by $Q_*(w)=S_*(w)I_2$, where $S_*\in\mathbf{A}$ 
is defined by $a(w)+b(w)S_*(w)=F_\gamma(p(w))$
for some non-constant $p\in\mathbf{RA}$ satisfying
$|p(w)|=1$ for all $w\in\mathbb T$,
and
$F_\gamma\in\mathbf A$ is defined in Lemma~\ref{lem1}.
\end{thm}

Theorem~\ref{thm0} corresponds to a special case of
Theorem~\ref{thm1}, with 
\[ a(w)=0.5w,\quad b(w)=w^2,\quad\gamma=\sqrt{2},\quad
p(w)=w.\]

A sketch of a
proof of Theorem~\ref{thm0} is given in the Appendix section below.

\section{Appendix}
The appendix contains
proof of the main results (Theorems 1 and 2), as well as that
of the (well known) statement of Lemma 1.

\subsection{Proof of Lemma~\ref{lem1}}
For $\theta\in(0,\pi)$ let
\[ \mathbb C_{\theta}=\left\{ re^{jt}:~ r>0,~|t|<\theta\right\},
\quad
\hat{\mathbb C}_{\theta}=\left\{ re^{jt}:~ r\ge0,~|t|\le\theta\right\}
\cup\{\infty\}\]
denote the "open angle $2\theta$ cone" in $\mathbb C$ and its
closure in $\mathbb C\cup\{\infty\}$. 
By definition, 
$F_\gamma=U_\alpha\circ R_\alpha\circ V$ 
is a composition of one power function
$R_\alpha:~\mathbb C_{\pi/2}\to\mathbb C_{\alpha}$
and
two M\"obius transformations 
$V:~\mathbb D\to\mathbb C_{\pi/2}$, 
$U_\alpha:~\mathbb C_{\alpha}\to\Omega_\gamma$
defined by
\[ V(w)=\frac{1+jw}{1-jw},\qquad
R_\alpha(s)=s^{2\alpha/\pi},\qquad
U_\alpha(y)=j\tan(\alpha)\cdot\frac{1-y}{1+y}~.\]
Since each function $V,U_\alpha,R_\alpha$ is a holomorphic
bijection, $F_\gamma$ is a holomorphic bijection, too.
Moreover, since $V,U_\alpha,R_\alpha$ have continuous
extensions $\hat V:~\mathbb D\cup\mathbb T\to
\hat{\mathbb C}_{\pi/2}$, 
$\hat U_\alpha:~\hat{\mathbb C}_{\alpha}\to\hat\Omega_\gamma$,
$\hat R_\alpha:~\hat{\mathbb C}_{\pi/2}\to\hat{\mathbb C}_{\alpha}$
(where $\hat\Omega_\gamma$ is the 
closure of $\Omega_\gamma$), $F_\gamma$ has a continuous 
extension $\hat F_\gamma:~\mathbb D\cup\mathbb T\to
\hat\Omega_\gamma$.
In addition, while the maps $V$ and $U_\alpha$ do not have real
symmetry, they satisfy conditions
\[  V(\bar w)=(\overline{V(w)})^{-1},\qquad
R_\alpha(1/\bar s)=(\overline{R_\alpha(s)})^{-1},\qquad
U_\alpha(1/\bar y)=\overline{U_\alpha(y)},\]
which proves that the total composition 
$F_\gamma=U_\alpha\circ R_\alpha\circ V$ 
is real symmetric. Finally, the expression for $\dot F_\gamma(0)$
follows from the observation that 
\[  \dot V(0)=2j,\quad V(0)=1,\quad \dot R_\alpha(1)=2\alpha/\pi,\quad
R_\alpha(1)=1,\quad\dot U_\alpha(1)=-0.5j\tan(\alpha).
\]

\subsection{Proof of Theorem~\ref{thm0}}
Since
\[ L_0(w)+L_1(w)Q(w)=
\begin{bmatrix}1&0.5w+w^2S_2(w)\\ 
0.5w+w^2S_1(w)& 1\end{bmatrix}
\quad\text{for}\quad
Q=\begin{bmatrix}S_1&0\\ 0& S_2\end{bmatrix},\]
the set of all transfer matrices $L_0+L_1Q$ 
with $Q\in \mathbf H^\infty_{2\times2}\cap\mathcal D$ can
be represented in the form
\[  \left\{L_0+L_1Q:~Q\in \mathbf H^\infty_{2\times2}\cap\mathcal D
\right\}=\left\{ H[G_1,G_2]:~G_1,G_2\in\mathbf X\right\},\]
where
\[  \mathbf X=
\left\{G\in\mathbf H^\infty:~G(0)=0,~\dot G(0)=0.5\right\},\qquad
H[G_1,G_2]=
\begin{bmatrix}1&G_2\\ G_1& 1\end{bmatrix}.\]
Theorem~\ref{thm0} claims that 
the infimum of $\|H[G_1,G_2]\|_\infty$ over 
$G_1,G_2\in\mathbf X\cap\mathbf{RA}$ equals 
$\sqrt2$, same as the
minimum of $\|H[G_1,G_2]\|_\infty$
over $G_1,G_2\in\mathbf X$, which in turn is achieved at
a unique pair $G_1=G_2=F_{\sqrt2}$.
The proof is presented in several steps.

\paragraph{Step 1.}
Note  that $\mathbf X$ is an affine subspace in
$\mathbf H^\infty$, i.e. 
$G=0.5(G_1+G_2)\in\mathbf X$ whenever 
$G_1,G_2\in\mathbf X$. 
Moreover, since
\[  \|H[G_2,G_1]\|_\infty=\|H[G_1,G_2]\|_\infty,\qquad
H[G,G]=0.5(H[G_1,G_2]+H[G_2,G_1]),\] 
convexity of the
H-Infinity norm function implies
$\|H[G_1,G_2]\|_\infty\ge\|H[G,G]\|_\infty$.
In other words, $\|H[G_1,G_2]\|_\infty\le\gamma$ for
$G_1,G_2\in\mathbf X$ implies $\|H[G,G]\|_\infty\le\gamma$
for $G=0.5(G_1+G_2)\in\mathbf X$.

\paragraph{Step 2.}
Since
\[  \sigma_{\max}\left(\begin{bmatrix}1&g\\ g&1\end{bmatrix}
\right)=\max\{|1-g|,|1+g|\}\qquad\forall~g\in\mathbb C,\]
the inequality
$\|H[G,G]\|_\infty\le\gamma$, where
$\gamma>1$ and $G\in\mathbf H^\infty$ is not constant
(note that {\sl all} $G\in\mathbf X$ are not constant) holds
if and only if
$G(w)\in\Omega_\gamma$
for all $w\in\mathbb D$.  

\paragraph{Step 3.}
Whenever $G\in\mathbf X$ is such that
$G(w)\in\Omega_\gamma$
for all $w\in\mathbb D$,
the composition 
$p=F_\gamma^{-1}\circ G$ satisfies conditions
\begin{equation}\label{eq42}
p\in\mathbf H^\infty,\quad \|p\|_\infty\le1,\quad
p(0)=0,
\end{equation}
\begin{equation}\label{eq43}
 \dot p(0)=\frac1{2\dot F_\gamma(0)}=
\frac{\pi\cos\alpha}{4\alpha\sin\alpha}\qquad
(0<\alpha<\pi/2,~\cos\alpha=1/\gamma)
\end{equation}
whenever $\|H[G,G]\|_\infty\le\gamma$.

\paragraph{Step 4.}
Since the Cauchy integral identity yields
\[   \dot p(0)=\frac1{2r\pi}\int_{-\pi}^\pi 
e^{-jt}p\left(re^{jt}\right)dt\]
for every $p\in\mathbf H^\infty$ and $r\in(0,1)$, it follows that
$|\dot p(0)|\le1$ whenever $p$ satisfies conditions
(\ref{eq42}), with equality $\dot p(0)=1$ possible
only when $p(w)\equiv w$.
Hence $\alpha\le\pi/4$ in (\ref{eq43}), 
i.e. $\gamma\ge\sqrt2$ whenever 
$\|H(G_1,G_2)\|_\infty\le\gamma$ for
$G_1,G_2\in\mathbf X$, with equality
$\|H(G_1,G_2)\|_\infty=\sqrt2$ possible only when
$0.5(G_1+G_2)=F_{\sqrt2}$. In particular,
$F_{\sqrt2}\in\mathbf X$, and 
$\|H(F_{\sqrt2},F_{\sqrt2})\|_\infty=\sqrt2$.

\paragraph{Step 5.}
As established at step 3, the functional $\|H[G_1,G_2]\|_\infty$
achieves its minimal value over $G_1,G_2\in\mathbf X$
when $G_1=G_2=F_{\sqrt2}$. To show that this is the {\sl only}
argument of minimum, let $G_1,G_2\in\mathbf X$ be
any pair satisfying $\|H[G_1,G_2]\|_\infty=\sqrt2$.
Then $G=0.5(G_1+G_2)=F_{\sqrt2}$. Let 
$T_1=H[G_1,G_2]$, $T_2=H[G_2,G_1]$,
$D=0.5(G_1-G_2)$.
Applying matrix identity
\[  M'_aM_a+M_d'M_d=\frac{M_1'M_1+M_2'M_2}2,\quad
\text{where}
\quad M_a=\frac{M_1+M_2}2,\quad
M_d=\frac{M_1-M_2}2\]
to
$M_1=T_1(w)$ and $M_2=T_2(w)$, with $w\in\mathbb D$,
in which case the diagonal elements of $M_i'M_i$ are not
larger than 2, the diagonal elements of
$M_a'M_a$ equal $1+|w|^2$, and 
the diagonal elements of  $M_d'M_d$ equal $|D(w)|^2$,
we conclude that $|D(w)|^2\le1-|w|^2$ which, due to the
maximum modulus principle, implies $D(w)\equiv0$,
i.e. $G_1=G_2=F_{\sqrt2}$.

\paragraph{Step 6.}
Finally,
to show that the maximal lower bound of $\|H[G_1,G_2]\|_\infty$
over $G_1,G_2\in\mathbf{RA}\cap\mathbf X$ equals $\sqrt2$,
note that $F_{\sqrt2}$ (as any other function from
class $\mathbf{A}$) can be approximated arbitrarily well
by polynomials, i.e. for every $\delta>0$ there exists a polynomial
$R_\delta\in\mathbf{RA}$ such that 
$\|R_\delta-F_{\sqrt2}\|_\infty<\delta$. Then $R_\delta(0)\to0$
and $\dot R_\delta(0)\to0.5$ as $\delta\to0$, hence 
$\|F_{\sqrt2}-G_\delta\|_\infty\to0$ and
$\|H[G_\delta,G_\delta]\|_\infty\to\sqrt2$ for
\[ G_\delta=\frac{R_\delta-R_\delta(0)}{2\dot R_\delta(0)}
\in\mathbf X\qquad(0<\delta<0.5,~\delta\to0).\]

\subsection{Proof of Theorem~\ref{thm1} (a sketch)}
We follow the lines of the proof of Theorem~\ref{thm0},
with some minor modifications.
We re-define
$\mathbf X=\left\{a+bS:~S\in\mathbf H^\infty\right\}$.
Theorem~\ref{thm1} claims that 
the infimum $\gamma$ of $\|H[G_1,G_2]\|_\infty$ over 
$G_1,G_2\in\mathbf X\cap\mathbf{RA}$ is always greater than 1, and
equals the
minimum of $\|H[G_1,G_2]\|_\infty$
over $G_1,G_2\in\mathbf X$, which in turn is achieved at
a unique pair $G_1=G_2=F_{\gamma}\circ p$, where $p\in\mathbf{RA}$
is not constant, and satisfies $|p(z)|=1$ for all $z\in\mathbb T$.

The reduction of the task of minimizing $\|H[G_1,G_2]\|_\infty$ over 
$G_1,G_2\in\mathbf X$ to the minimization task
\begin{equation}\label{eq47}
r\to\min,\qquad\text{subject to}\qquad
G\in\mathbf X,\qquad G(w)\in\Omega_r\quad\forall~w\in\mathbb D
\end{equation}
is done the same way as
 in the proof of Theorem~\ref{thm0}. Since
the sets $\{G\in\mathbf H^\infty:~\|G\|_\infty\le R\}$ are compact
in the topology of uniform convergence on all compact subsets of
$\mathbb D$, and since the function 
$G\to\|G\|_\infty$ is lower semi-continuous
in this topology, there exists an
optimal  $G=G_*\in\mathbf X$. Then, for $\gamma=\min~r$,
the function
$p=F_\gamma^{-1}\circ G_*\in\mathbf H^\infty$ satisfies
$|p(w)|<1$ for all $w\in\mathbb D$, i.e. $\|p\|_\infty\le1$. 
Moreover, there exist no $S\in\mathbf H^\infty$
such that $\|p+bS\|_\infty<1$, because otherwise
$G=F_\gamma\circ(p+bS)\in\mathbf X$ would satisfy
the constraints in (\ref{eq47})
for some $r<\gamma$.
According to the standard theory of frequency-domain
H-Infinity optimization
(see, for example, \cite{Francis87}), $p$ is a 
rational function of
order smaller than the number of roots of $b$ in $\mathbb D$
(counting multiplicity), which satisfies the condition
$|p(z)|=1$ for all $z\in\mathbb T$, and is unique in the sense that
$\|p+bS\|_\infty>1$ for all non-zero $S\in\mathbf H^\infty$.
Since every other minimizer $G$ 
in (\ref{eq47}) 
will satisfy
$S=(F_\gamma^{-1}\circ G-p)/b\in\mathbf H^\infty$, this confirms 
$G=G_*=F_\gamma^{-1}\circ p$ as the unique minimizer 
in (\ref{eq47}).

\bibliography{notrational}{}
\bibliographystyle{plain}

\end{document}